\documentclass[11pt]{amsart}
\usepackage{amsmath,amsthm,epsfig,amssymb}
\usepackage{xypic}
\input xy
\title{Counting all regular octahedrons in $\{0,1,...,n\}^3$}
\author{Eugen J. Ionascu}
\curraddr{Department of Mathematics\\ Columbus State University\\4225 University Avenue\\
Columbus, GA 31907\\
Honorific Member of the Romanian Institute of Mathematics ``Simion
Stoilow" } \email{Ionascu\_Eugen@colstate.edu}

\date{July $9^{th}$, 2010}
\begin{document}
\def\sms{\small\scshape}
\baselineskip18pt
\newtheorem{theorem}{\hspace{\parindent}
T{\scriptsize HEOREM}}[section]
\newtheorem{proposition}[theorem]
{\hspace{\parindent }P{\scriptsize ROPOSITION}}
\newtheorem{corollary}[theorem]
{\hspace{\parindent }C{\scriptsize OROLLARY}}
\newtheorem{lemma}[theorem]
{\hspace{\parindent }L{\scriptsize EMMA}}
\newtheorem{definition}[theorem]
{\hspace{\parindent }D{\scriptsize EFINITION}}
\newtheorem{problem}[theorem]
{\hspace{\parindent }P{\scriptsize ROBLEM}}
\newtheorem{conjecture}[theorem]
{\hspace{\parindent }C{\scriptsize ONJECTURE}}
\newtheorem{example}[theorem]
{\hspace{\parindent }E{\scriptsize XAMPLE}}
\newtheorem{remark}[theorem]
{\hspace{\parindent }R{\scriptsize EMARK}}
\renewcommand{\thetheorem}{\arabic{section}.\arabic{theorem}}
\renewcommand{\theenumi}{(\roman{enumi})}
\renewcommand{\labelenumi}{\theenumi}
\newcommand{\Q}{{\mathbb Q}}
\newcommand{\Z}{{\mathbb Z}}
\newcommand{\N}{{\mathbb N}}
\newcommand{\C}{{\mathbb C}}
\newcommand{\R}{{\mathbb R}}
\newcommand{\F}{{\mathbb F}}
\newcommand{\K}{{\mathbb K}}
\newcommand{\D}{{\mathbb D}}
\def\phi{\varphi}
\def\ra{\rightarrow}
\def\sd{\bigtriangledown}
\def\ac{\mathaccent94}
\def\wi{\sim}
\def\wt{\widetilde}
\def\bb#1{{\Bbb#1}}
\def\bs{\backslash}
\def\cal{\mathcal}
\def\ca#1{{\cal#1}}
\def\Bbb#1{\bf#1}
\def\blacksquare{{\ \vrule height7pt width7pt depth0pt}}
\def\bsq{\blacksquare}
\def\proof{\hspace{\parindent}{P{\scriptsize ROOF}}}
\def\pofthe{P{\scriptsize ROOF OF}
T{\scriptsize HEOREM}\  }
\def\pofle{\hspace{\parindent}P{\scriptsize ROOF OF}
L{\scriptsize EMMA}\  }
\def\pofcor{\hspace{\parindent}P{\scriptsize ROOF OF}
C{\scriptsize ROLLARY}\  }
\def\pofpro{\hspace{\parindent}P{\scriptsize ROOF OF}
P{\scriptsize ROPOSITION}\  }
\def\n{\noindent}
\def\wh{\widehat}
\def\eproof{$\hfill\bsq$\par}
\def\ds{\displaystyle}
\def\du{\overset{\text {\bf .}}{\cup}}
\def\Du{\overset{\text {\bf .}}{\bigcup}}
\def\b{$\blacklozenge$}

\def\eqtr{{\cal E}{\cal T}(\Z) }
\def\eproofi{\bsq}

\begin{abstract} In this paper we describe a procedure for calculating the number
of regular octahedrons that have vertices with coordinates in the
set $\{0,1,...,n\}$. As a result, we introduce a new sequence in
{\it The Online Encyclopedia of Integer Sequences} (A178797) and
list the first one hundred terms of it. We adapt the method
appeared in \cite{ejicregt} which was used to find the number of
regular tetrahedra with coordinates of their vertices in
$\{0,1,...,n\}$. The idea of this calculation is based on the
theoretical results obtained in \cite{ejiam}. A new fact proved
here helps increasing the speed of all the programs used before.
The procedure is put together in a series of commands written for
Maple.
\end{abstract} \maketitle
\section{INTRODUCTION}

In this note we complete the work begun in the sequence of papers
\cite{rceji}, \cite{eji}-\cite{ejiam} about equilateral triangles,
regular tetrahedra, cubes, and regular octahedrons all with
vertices having integer coordinates. Very often we will refer to
this property by saying that the various objects are in $\mathbb
Z^3$. Strictly speaking these geometric objects are defined as
being more than the set of their vertices that determines them,
but for us here, these are just the vertices. So, for instance, an
equilateral triangle is going to be a set of three points in
$\mathbb Z^3$ for which the Euclidean distances between every two
of these points are the same. The main purpose of the paper is to
take a close look at the regular octahedrons in $\mathbb Z^3$. The
simplest example of a regular octahedron with integer coordinates
for its vertices is
$$OC_1:=\{[0, 1, 1], [1, 0, 1], [1, 1, 0], [1, 1, 2], [1, 2, 1], [2,
1, 1]\},$$ \n which can be obtained from the usual unit cube in
$\mathbb R^3$, multiplying the vertices by a factor of two and
then taking the coordinates of the centers of the new faces. It
turns out that this procedure gives all such octahedrons as shown
in \cite{ejiam}:

\begin{center}\label{fig1}
$\underset{\small \ Figure\ 1\ (a) : \ OC_1\
octahedron}{\epsfig{file=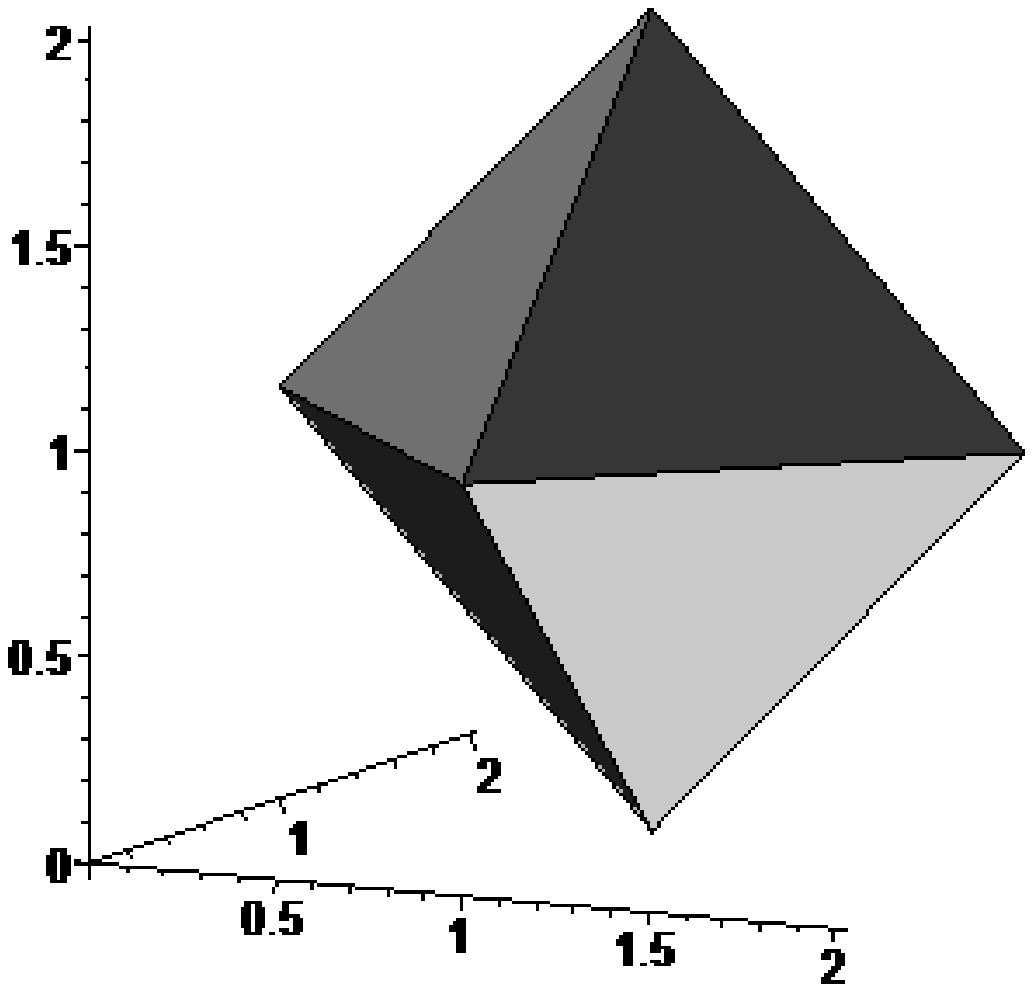,height=2in,width=2in}}$
$\underset{\small \ Figure\ 1\ (b) : \ Regular\ octahedron vs
cube}{\epsfig{file=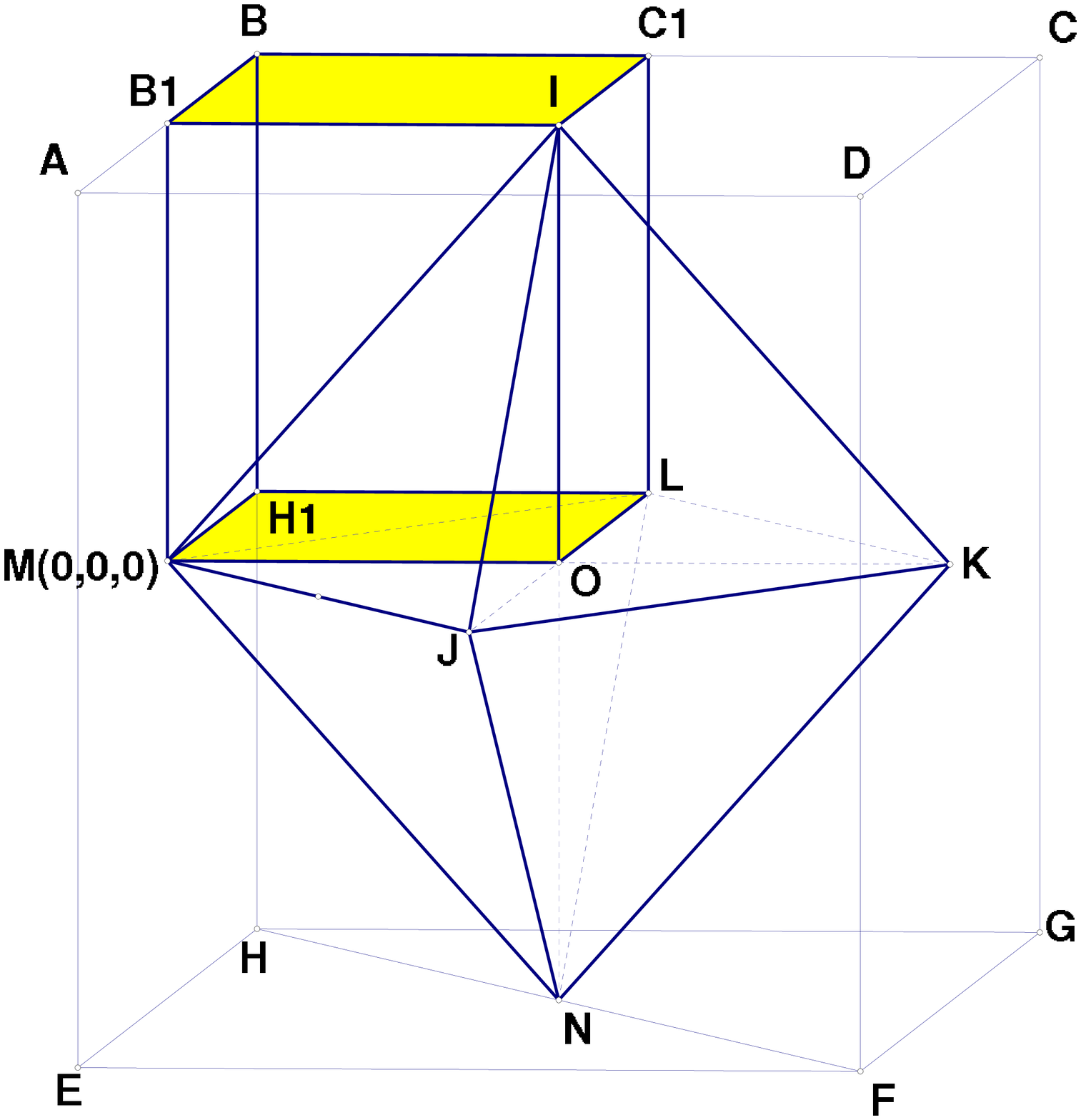,height=2in,width=2in}}$
\end{center}

\begin{theorem}\label{oldtheorem} Every regular octahedron in $\mathbb Z^3$ is the
dual of a cube that can be obtained (up to a translation with a
vector with integer coordinates) by doubling a cube in $\mathbb
Z^3$.
\end{theorem}

Referring to Figure~\ref{fig1} (b), we showed that if the regular
octahedron $IJKLMN$ is in  $\mathbb Z^3$, then so is the cube
$BB_1C_1IH_1LOM$ and vice versa. This defines a one-to-one
correspondence between the classes of cubes (invariant under
integer translations) and the classes regular octahedra (invariant
under integer translations) in $\mathbb Z^3$. In \cite{ejio} we
determined a sequence of irreducible cubes, one from each of the
classes of cubes invariant under integer translations and cube
symmetries. For each one of these cubes we can construct as before
a regular octahedron, obtaining this way a sequence of irreducible
regular octahedrons.

\section{Some new ingredients and other theoretical facts}

In \cite{ejio}, we improved and adapted the earlier code for
counting all cubes with vertices in $\{0,1,...,n\}^3$ and extended
the sequence A098928.  In this paper, the usual techniques and
ideas are going to be the same except some counting procedure that
is very efficient in comparison to what we had before. We are
going to treat this in the general case so, let us suppose that
these objects can be either equilateral triangles, regular
tetrahedrons, cubes or regular octahedrons with vertices in
$\mathbb Z^3$. For such an object, say $\cal O$, we can translate
it, within $\mathbb Z^3$, to $\cal O'$ that is in the positive
octant and in such way each plane of coordinates contains at least
one vertex of $\cal O'$. Let us denote by $\alpha_0$ the number of
objects in $C_m$ obtained by applying to $\cal O'$  all 48
possible symmetries of the cube $C_m$. These symmetries are
generated in the following way: first we have symmetries with
respect to the middle planes and compositions, for example

$$(x,y,z)\to (m-x,y,z), \ (x,y,z)\to (m-x,m-y,z),\  (x,y,z)\to (m-x,m-y,m-z),$$

\n in a total of eight including the identity, then each one of
these is coupled with one of the six permutations of the variables
($\cal S_6$). These transformations form a group isomorphic with
the group of all $3\times 3$ orthogonal matrices having entries
$\pm 1$ and it is also known as the group of symmetries of a cube
or of a regular octahedron. It is isomorphic to $S_4\times \mathbb
Z_2$. We are going to denote this group by ${\cal S}_{cube}$
although it is usually known under the name of (extended)
octahedral group and denoted simply by $O_h$.

If we think of $\alpha_0$ as the cardinality of

$$Orbit(\cal O'):=\{ s({\cal O'}) |s\in {\cal
S}_{cube}\},$$

\n which is, by the first theorem of isomorphism of groups, the
same as the cardinality of the group factor ${\cal S}_{cube}/{\cal
G}$, where ${\cal G}$ is the subgroup of ${\cal S}_{cube}$ of
those symmetries that leave $\cal O'$ invariant. The structure of
subgroups of ${\cal S}_{cube}$ is known and for each divisor of 48
there is a subgroup of that order. Hence, we expect $\alpha_0$ to
be in the set $\{1, 2, 3, 4, 6, 8, 12, 16, 24, 48\}$ and most of
the time to be $48$ since an arbitrary object $\cal O'$ in $C_m$
is unlikely to be invariant under any of the symmetries of ${\cal
S}_{cube}$.

Then, we denote by $\alpha$, the cardinality of the set of all the
objects counted in $\alpha_0$ and their (all possible) integer
translations that leave the resulting objects in $C_m$. Also, we
denote by $\beta$ the objects counted in $\alpha$ which are in
$\{0,1,...,m\}^2\times \{0,1,...,m-1\}$. Finally, let us denote by
$\gamma$ the number of objects counted in $\beta$ which are in
$\{0,1,...,m\}\times \{0,1,...,m-1\}^2$. Then, we found a formula
that gives the number of objects obtained from $\cal O$, under all
symmetries and translation that leaves the resulting object in
$\{0,1,...,k\}^3$, $k\ge m$.

This fact has been essentially proved in Theorem 2.2 in
\cite{ejic}. The formula that gives this number is

\begin{equation}\label{oldformula}
N({\cal
O},k)=(k-m+1)^3\alpha-3(k-m)(k-m+1)^2\beta+3(k-m+1)(k-m)^2\gamma.
\end{equation}

Let us suppose that the object $\cal O$ can be squeezed within a
box of dimensions $m\times n\times p$ ($m\ge n\ge p$), i.e. up to
symmetries and translations, $\cal O$ can be transformed to $\cal
O'$ fitting snugly into
$$B_{m,n,p}:=\{0,1,...,m\}\times \{0,1,...,n\}\times \{0,1,...,p\}.$$

We can similarly consider all eight reflections compatible with
the box $B_{m,n,p}$ of the form

$$(x,y,z)\to (m-x,y,z), \ (x,y,z)\to (m-x,n-y,z),\  (x,y,z)\to (m-x,n-y,p-z), etc. $$

Let us denote the group of these transformations by $\cal S_b$. We
notice that each one of these transformation leaves the object
$\cal O'$ inside the box $B_{m,n,p}$. From case to case, depending
of what the values $m$, $n$ and $p$ are, we may have the result of
some or all of the permutation transformations applied to $\cal
O'$ still in $B_{m,n,p}$. Hence, we will denote by $\omega(\cal
O)$ the cardinality of the set

$$BoxOrbit(\cal O'):=\{ [s_1\circ s_2] ({\cal O'}) \in B_{m,n,p}|s_1\in {\cal
S}_b, s_2\in {\cal S}_6\}.$$

Let us look at an example. Suppose $\cal O$ (equal with $\cal O'$)
is the equilateral triangle given by its vertices:  $$\{[0, 2, 2],
[5, 7, 0], [7, 0, 1]\}.$$

We observe that ${\cal O}\in B_{7,7,2}$. Then one can check that
$BoxOrbit(\cal O)$ is the collection of eight triangles

$${\small \begin{array}{c} {\cal O}, \{[0, 0, 1], [2, 7, 0], [7, 2, 2]\},
\{[0, 0, 1], [2, 7, 2], [7, 2, 0]\}, \{[0, 2, 0], [5, 7, 2], [7,
0, 1]\},\\ \\  \{[0, 5, 0], [5, 0, 2], [7, 7, 1]\}, \{[0, 5, 2],
[5, 0, 0], [7, 7, 1]\}, \{[0, 7, 1], [2, 0, 0], [7, 5, 2]\}, \{[0,
7, 1], [2, 0, 2], [7, 5, 0]\},
\end{array}}$$

\n so $\omega(\cal O)=8$. It turns out that $\alpha_0(\cal O)=48$,
$\alpha(\cal O)=144$, $\beta(\cal O)=40$ and $\gamma(\cal O)=0$.
Formula (\ref{oldformula}) becomes

$$N({\cal O},k)=24(k-1)(k-6)^2,\ \ k\ge 7. $$

It turns out the this factorization is not accidental and the
following alternative to (\ref{oldformula}) is true.

\begin{theorem}\label{coj}
Given $\cal O$, one of the objects mentioned before, and
$B_{m,n,p}$ the smallest box containing a translation of $\cal O$
($m\ge n\ge p$), we let $u=m-n$, $v=n-p$, and
$$\Delta=\omega({\cal O})(k-m+1)(k-n+1)(k-p+1).$$
Then the number of distinct objects in the cube $B_{k,k,k}$ ($k\ge
m$), obtained from $\cal O$ by all possible integer translations
and symmetries is equal to
\begin{equation}\label{newformula}
N({\cal O},k)=\begin{cases} \Delta \ if \ u=v=0,\\ \\
3\Delta\ if \ u \ or\ v \ is \ 0,\\ \\
6\Delta \ if \ u\  and \ v >0.
\end{cases}
\end{equation}
\end{theorem}

\proof.\  The case $u=v=0$ implies $\omega(\cal O)=\alpha_0(\cal
O)=\alpha(\cal O)$ and $\beta(\cal O)=\gamma(\cal O)=0$ because
there is no room to shift the orbit $Orbit(\cal O')$ inside of
$B_{m,m,m}$. The formula follows from (\ref{oldformula}).

Let us look into the case $u>0$ and $v>0$.  We begin by observing
that each integer translation of the box $B_{m,n,p}$ in all
possible ways inside $B_{k,k,k}$ will give $\omega(\cal O)$ more
copies of $\cal O$. There is no overlap between these copies
because neither one of them can be inside of two distinct
translations of $B_{m,n,p}$. This is due to the minimality of $m$,
$n$ and $p$. We get $\Delta$ such copies by counting all possible
translations. Since $m$, $n$ and $p$ are all distinct, the box
$B_{m,n,p}$ can be positioned first with the biggest of its
dimensions along one of the directions given by the axis of
coordinates, that is three different ways, and for each such
position the next largest dimension can be positioned along the
two remaining directions. The minimality of $m$, $n$ and $p$ makes
the six  different situations generate distinct objects. This
explains the factor of six that appears in (\ref{newformula}) for
this situation.

In the last case, the box $B_{m,n,p}$ has two of its dimensions
the same, so there are only three possibilities to arrange the box
before one translates it. To see that we get all possible
translates and symmetries of $\cal O$ by this counting, we can
start with one copy $\cal O'$. Construct the minimum box around
it. In terms of its position and dimensions,  we know in what of
the six or three cases we are. We transform it into the standard
standard position, $B_{m,n,p}$, and look at the corresponding
object, $\cal O''$. The transformations involved form a group of
transformations generated by the permutations of the coordinates,
the reflections into the axes and integer translations. Every
transformation in this group, say $g=\tau\circ \sigma\circ \pi$
with $\pi$ a permutation, $\sigma$ a reflection or a composition
of reflections and $\tau$ a translation, which satisfies $g(\cal
O)=\cal O''$ determines a representation $(s_1\circ s_2)(\cal
O)=\cal O''$ with $s_1\in {\cal S}_b$, $s_2\in {\cal S}_6$ as in
the definition of $\omega(\cal O)$. This can be done by taking
$s_2=\pi$ and $s_1=\tau\circ \sigma$. This is true again because
of the minimality of the box $B_{m,n,p}$, i.e. there is only one
integer translation that takes a reflected box $B'_{m,n,p}$ into
$B_{m,n,p}$.\eproof

This new way of counting is more efficient from a computational
point of view because $\omega$ is simply no bigger than 48, as
opposed to the previous situation when $\alpha$, $\beta$ and
$\gamma$ could turn out to be big numbers and so the number of
iterations for computing them would be also large. Roughly
speaking, this counting factors out fast the problem with the
integer translations.
%
%
%
%
%
%

As an example, let us consider $$OC_2=\{[0, 0, 1], [0, 3, 4], [1,
4, 0], [3, 0, 4], [4, 1, 0], [4, 4, 3]\}.$$

The minimal box here is $B_{4,4,4}$ and after rotating $OC_2$ in
all possible ways (Figure~\ref{fig2} (b)) we get $\omega(OC_2)=4$.

\begin{center}\label{fig2}
$\underset{\small \ Figure\ 2\ (a) : \ OC_2\
octahedron}{\epsfig{file=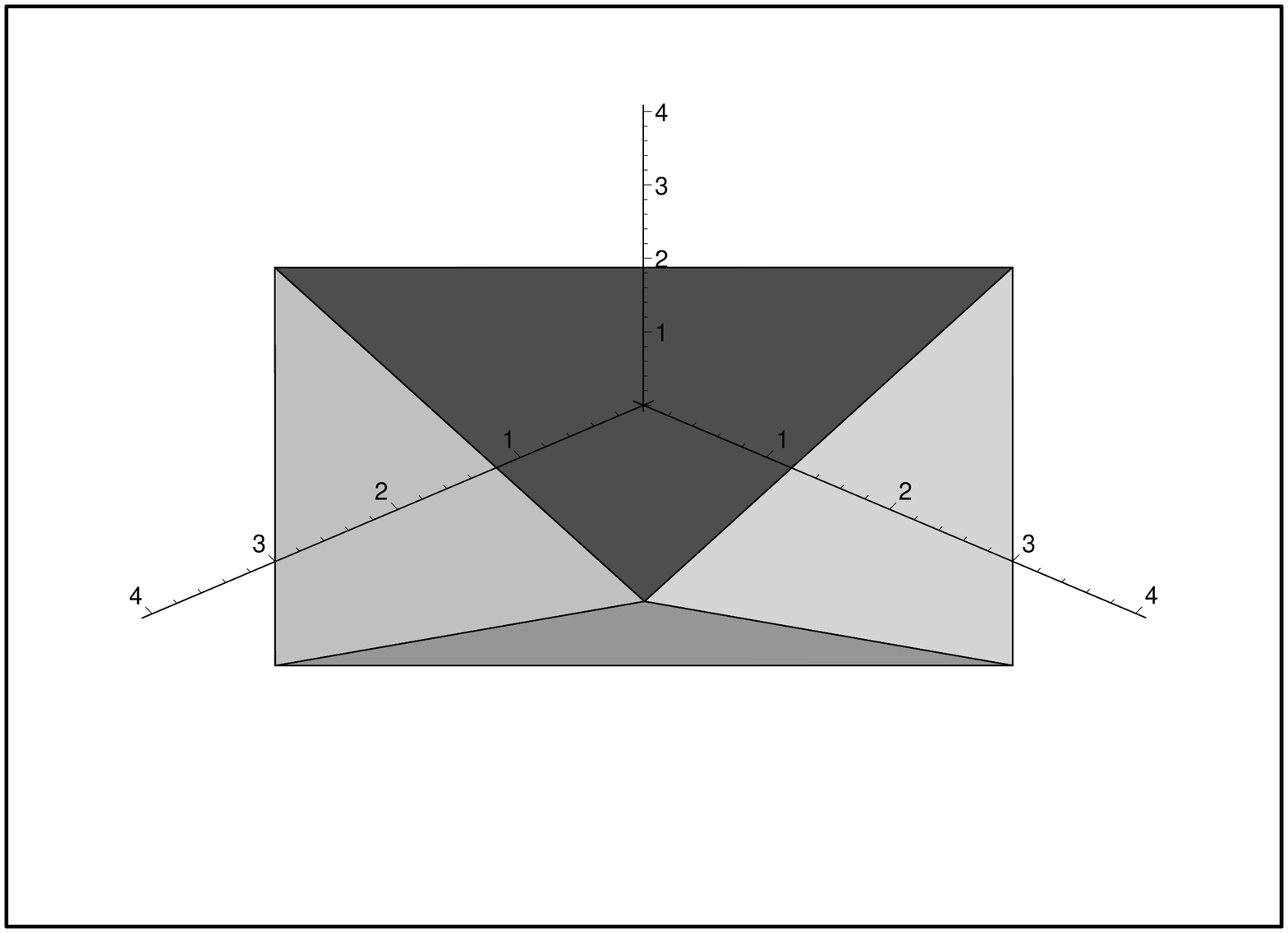,height=2in,width=2in,angle=-90}}$\
\  $\underset{\small \ Figure\ 2\ (b) : Four\ octahedrons\ in\
the\ box\
}{\epsfig{file=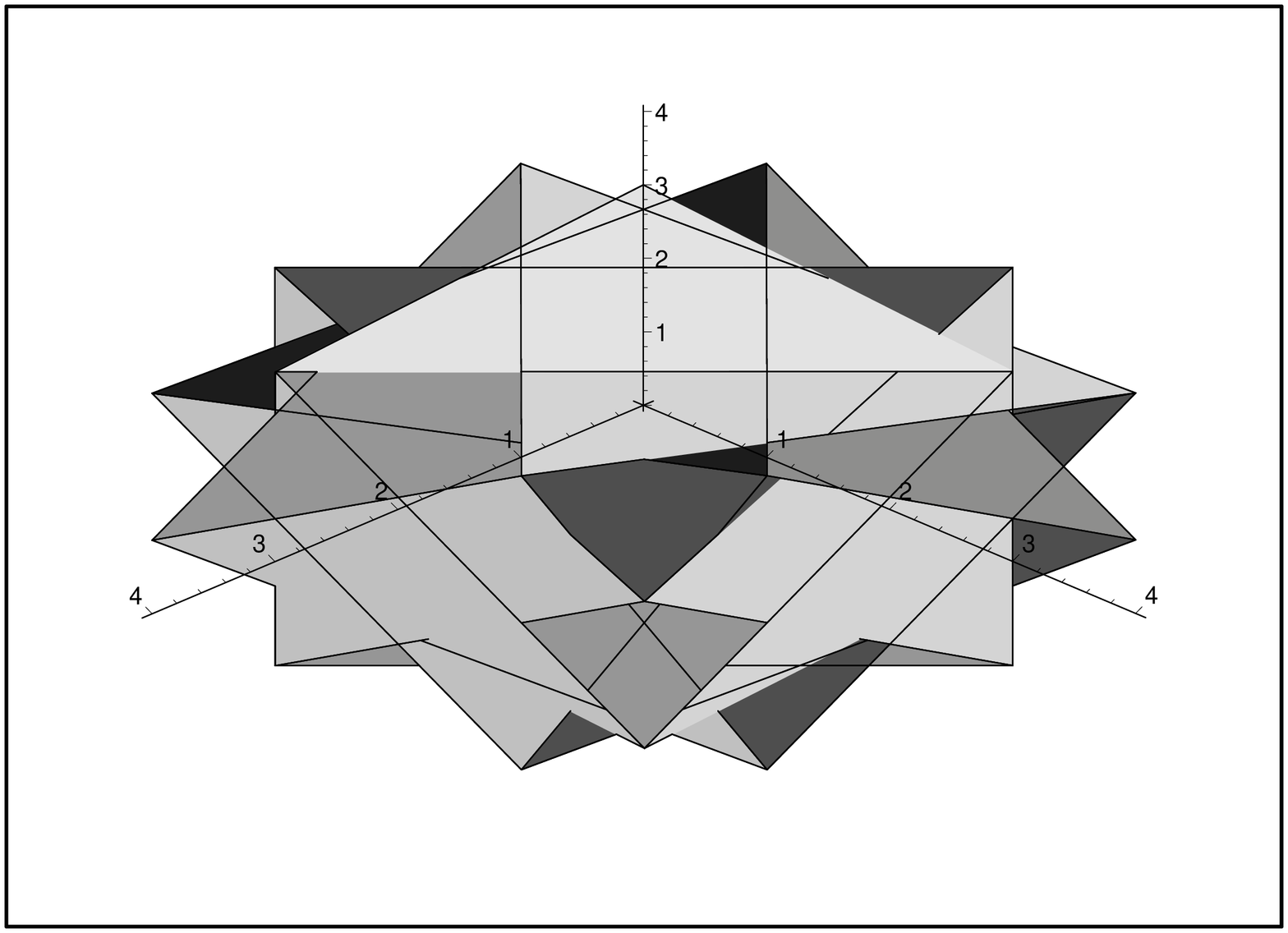,height=2in,width=2in,angle=-90}}$
\end{center}


The idea of calculations is basically the same as in \cite{ejio},
in which we have constructed a list of irreducible cubes that are
used to generate all the other cubes in $B_{k,k,k}$. Here, we are
using Theorem~\ref{oldtheorem}, to construct a similar list of
irreducible regular hexahedrons. As expected, an {\it irreducible}
regular hexahedron is one whose coordinates cannot be obtained
from a strictly smaller one with vertices in $\mathbb Z^3$ by
integer dilations and translations. One simple consequence of
Theorem~\ref{oldtheorem} is the next corollary.

\begin{corollary}\label{sideofoctahedrons}
The sides of every irreducible regular octahedron in $\mathbb Z^3$
are of the form $$(2k-1)\sqrt{2}$$\n with $k\in \mathbb  N$.
\end{corollary}

We used the same way of finding the parameterizations of the
equilateral triangles as in \cite{ejio}:

\begin{center}\label{fig3}
$\underset{\small \ Figure\ 3\ (a) : \ Vectors \ \eta\ and \ \zeta
}{\epsfig{file=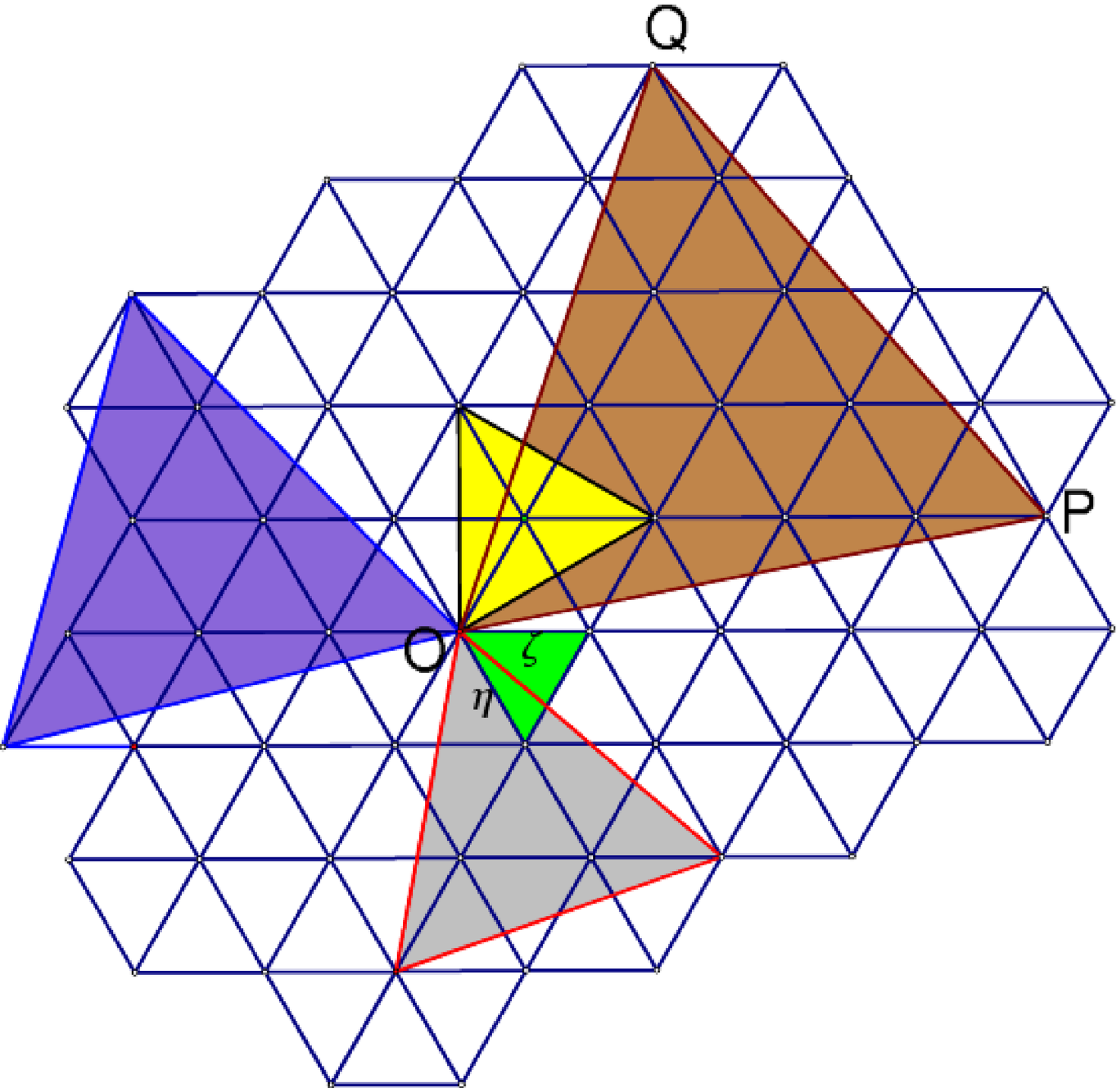,height=2in,width=2in}}$\ \
$\underset{\small \ Figure\ 3\ (b) : Minimality \ of \
parametrization
}{\epsfig{file=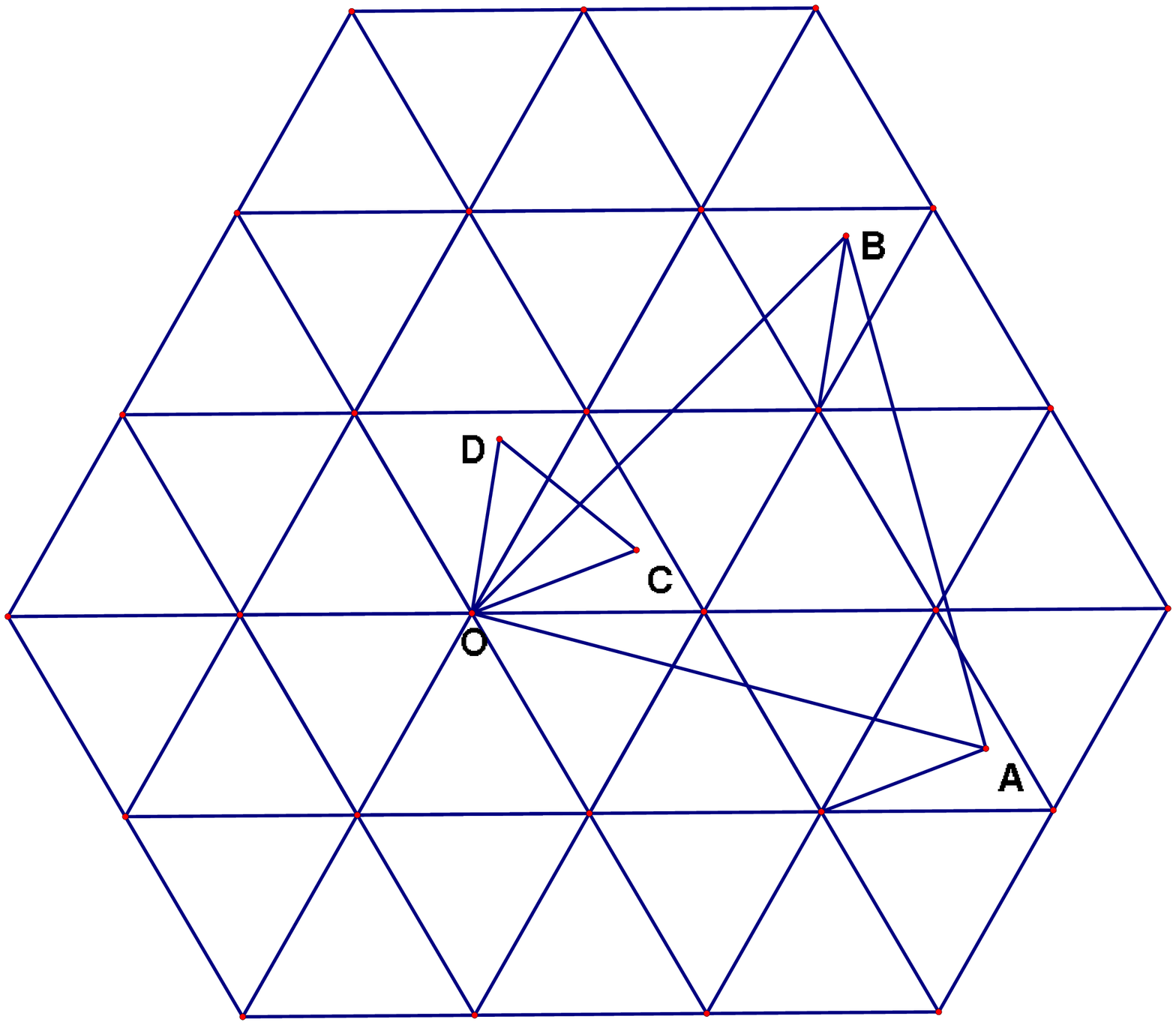,height=1.7in,width=2in}}$
\end{center}

\begin{equation}\label{vectorid}
\overrightarrow{OP}=m\overrightarrow{\zeta}-n\overrightarrow{\eta},\
\
\overrightarrow{OQ}=n\overrightarrow{\zeta}-(n-m)\overrightarrow{\eta},
\ { \rm with} \ \overrightarrow{\zeta}=(\zeta_1,\zeta_1,\zeta_2),
\overrightarrow{\eta}=(\eta_1,\eta_2,\eta_3),
\end{equation}

\begin{equation}\label{paramtwo}
\begin{array}{l}
\begin{cases}
\ds \zeta_1=-\frac{rac+dbs}{q}, \\ \\
\ds \zeta_2=\frac{das-bcr}{q},\\ \\
\ds \zeta_3=r,
\end{cases}
,\ \
\begin{cases}
\ds \eta_1=-\frac{db(s-3r)+ac(r+s)}{2q},\\ \\
\ds \eta_2=\frac{da(s-3r)-bc(r+s)}{2q},\\ \\
\ds \eta_3=\frac{r+s}{2},
\end{cases}
\end{array}
\end{equation}
where $q=a^2+b^2$ and $(r,s)$ is a suitable solution of
$2q=s^2+3r^2$ that makes all the numbers in {\rm (\ref{paramtwo})}
integers. The sides-lengths of $\triangle OPQ$ are equal to
$d\sqrt{2(m^2-mn+n^2)}$ and the triangle can be completed to a
regular tetrahedron (in space) with integer coordinates if and
only if $k^2=m^2-mn+n^2$ for some $k\in \mathbb Z$. Related to
this fact we have the following proposition.

\begin{proposition}\label{numberofro} For a prime $p>3$, the
number of irreducible regular octahedrons in $\mathbb Z^3$ (up to
translations and symmetries) having side lengths equal to
$p\sqrt{2}$, is at most $\pi\epsilon(p)+1$, where
\begin{equation}\label{numberofrepr}
\pi\epsilon(p)=\frac{\Lambda(p)+24\Gamma_2(3p^2)}{48},
\end{equation}
with $\Lambda$ and $\Gamma$ defined as in \cite{ejicregt}.
\end{proposition}

The exact number of such objects is yet a big mystery to us.

\section{The Maple Code}

We wrote the code using the Maple Software and so we took
advantage of the build in commands available for a various number
of functions. The beginning is pretty standard:

\begin{verbatim}
> restart:with(numtheory):with(plots):
\end{verbatim}

{\bf Step 1.}  The next three procedures calculate all possible
values of $k$, then the parameters $m$, $n$, and finally the
normal vector $(a,b,c)$. The values of $k$ are determined by using
a characterization theorem for the quadratic form involved here,
i.e. all values of $k$ less than $n$ such that $k$ is of the form
a product of primes of the form $3\ell+1$, $\ell \ge 2$ or $k=1$.

\begin{verbatim}
 > kvalues:=proc(n)
 local i,j,k,L,a,p,q,r,m,mm;
 L:={};mm:=floor((n+1)/2);
 for i from 2 to mm do
 a:=ifactors(2*i-1);
 k:=nops(a[2]);r:=0;
     for j from 1 to k do
       m:=a[2][j][1]; p:=m mod 3;
       if m=3 then r:=1 fi;
       if r=0 and p=2  then r:=1 fi;
     od;
 if r=0 then L:=L union {2*i-1};fi;
 od; L:=L union {1}; L:=convert(L,list);
 end:
\end{verbatim}

\n For example, $kvalues(100)=[1, 7, 13, 19, 31, 37, 43, 7^2, 61,
67, 73, 79, 7(13), 97]$.

\n  Next, we are interested in the solutions $(m,n)$  of the
equation $k^2=m^2-mn+n^2$, which are primitive in the sense that
$gcd(m,n)=1$, $m>0$, $n>0$, and $2m<n$. We apply this procedure to
only those $k$'s which are the output of $kvalues$.

\begin{verbatim}
> listofmn:=proc(k)
local a,b,i,nx,x,m,n,L;
x:=[isolve(k^2=m^2-m*n+n^2)];
 nx:=nops(x); L:={};
  for i from 1 to nx do
    if lhs(x[i][1])=m then a:=rhs(x[i][1]); b:=rhs(x[i][2]);
     else  b:=rhs(x[i][1]); a:=rhs(x[i][2]); fi;
    if  gcd(a,b)=1 and a>=0 and b>0 and 2*a<b then L:=L union {[a,b]};fi;
  od;L;
end:
\end{verbatim}

\n  A prime of the form $3\ell+1$ has a unique primitive
decomposition as described, for example,
$listofmn(79)=\{[40,91]\}$ since $79^2=40-40(91)+91^2$ and $79$ is
a prime. In general, the number of solutions is equal to
$2^{\omega(k)-1}$ where $\omega(k)$ is the number of distinct
factors of $k$ which are of the form $3\ell+1$. This choice of $m$
and $n$ helps identify the irreducible octahedrons: all the other
solutions $(m,n)$ of the equation $k^2=m^2-mn+n^2$ lead to the
same regular octahedron. We included the case $a=0$ to obtain the
output $\{[0,1]\}$ for $listofmn(1)$ which is necessary later on.

\n In the next procedure, we find only the solutions $(a,b,c)$ of
$a^2+b^2+c^2=3d^2$ that satisfy $gcd(a,b,c)=1$ and $0<a\le b\le
c$. The number of such solutions can be calculated precisely in
terms of the prime factorization of $d$ (see \cite{ejicregt}).

\begin{verbatim}
> abcsol:=proc(d)
local i,j,k,m,u,x,y,sol,cd;sol:={};
for i from 1 to d do
    u:=[isolve(3*d^2-i^2=x^2+y^2)];k:=nops(u);
    for j from 1 to k do
       if rhs(u[j][1])>=i and rhs(u[j][2])>=i then
         cd:=gcd(gcd(i,rhs(u[j][1])),rhs(u[j][2]));
       if cd=1 then sol:=sol union {sort([i,rhs(u[j][1]),rhs(u[j][2])])};fi;fi;
od; od;convert(sol,list); end:
\end{verbatim}

For $d=17$, for instance, the procedure {\it abcsol} gives the
four solutions: [1, 5, 29], [7, 17, 23], [11, 11, 25] and [13, 13,
23]. We observe that two of them have the property that $a=b$, in
which case, we know (see \cite{ejio}) that the formulae
(\ref{paramtwo}) simplify quite a bit.

\vspace{0.2in}

{\bf Step 2.}   The next seven procedures implement the new way of
finding $r$ and $s$ which appear in the parametrization
(\ref{paramtwo}). It is followed by the calculation of the of the
parametrization of equilateral triangles (\ref{paramtwo}) and only
one regular octahedron is constructed based on
Theorem~\ref{oldtheorem}. This octahedron is then translated
minimally into the positive octant and the minimal cube containing
it is computed.

For a prime  $p$ of the form $6\ell+1$, there exists an unique
decomposition $p=x^2+3y^2$ which is calculated next.

\begin{verbatim}
> uniquedecomposition:=proc(p)
local out,s,out1,out2;
if p=2 then out:=[1,1]; fi;
  if p>2 then s:=isolve(p=x^2+3*y^2);
   out1:=abs(rhs(s[1][1]));out2:=abs(rhs(s[1][2]));
    if out1^2+3*out2^2=p then out:=[out1,out2]; else out:=[out2,out1];fi;
  fi; out; end:
\end{verbatim}

For $p=2$ this procedure has the needed output $[1,1]$. As an
example,  $uniquedecomposition(2011)=[44,5]$ since $2011$ is a
prime and $2011=44^2+3(5^2)$. The next procedure is calculating
the factorization in $\mathbb Z[\sqrt{3} i]$ of a number of the
form $u+v\sqrt{3}i$.

\begin{verbatim}
> factoroverEisensteinintegers:=proc(u,v)
 local i,N,M,k,a,x,f1,f2,g,y,y1,y2,L,NN,MM,uu,vv;
 a:=sqrt(3)*I;NN:=gcd(u,v);uu:=u/NN;vv:=v/NN;
 N:=uu^2+3*vv^2;x:=uu+vv*a;M:=ifactors(N);k:=nops(M[2]);
 for i from 1 to k do
    f1:=M[2][i][1]; f2:=M[2][i][2];
      if f1>2 then
         g:=uniquedecomposition(f1);
         else  g:=[1,1]; f2:=1;
      fi;
    y:=expand(rationalize(x/(g[1]+a*g[2])));
    y1:=Re(y);y2:=type(y1,integer);
      if y2=true then L[i]:=[g[1]+a*g[2],f2]; else
        L[i]:=[g[1]-a*g[2],f2];
      fi;
od;[NN,seq(L[i],i=1..k),expand(NN*product(L[ii][1]^L[ii][2],ii=1..k))];end:
\end{verbatim}
As a simple example here, the following decomposition is obtained
for $u=13$ and $17$:

$$(1+I\sqrt{3})(2+I\sqrt{3})(5-2I\sqrt{3})=13+17I\sqrt{3}.$$
Perhaps one word of caution is necessary at this point. The
decomposition in general is not unique in the usual sense  since
$4=2(2)=(1+\sqrt{3})(1-\sqrt{3})$. As in the example shown, we go
for the second representation when something like this happens. We
are interested in this decomposition because it provides suitable
values for $r$ and $s$ that we find next. This turns out to be the
greatest common divisor between $A+I\sqrt{3}B$ and $2q$ with
$A=ac$, $B=bd$ and $q=a^2+b^2$ (see \cite{ejio}).

\begin{verbatim}
> findgcd:=proc(A,B,q)
local a,i,j,f,f1,f2,common,m,qq,fac,nfac,s,rs,P; a:=sqrt(3)*I;
f:=factoroverEisensteinintegers(A,B);m:=nops(f)-1;
common:=gcd(f[1],2*q); qq:=2*q/common^2; P:=common;
fac:=ifactors(qq);nfac:=nops(fac[2]);
   for i from 1 to nfac do
     f1:=fac[2][i][1];f2:=fac[2][i][2];
        if f1=2 then f2:=1;fi;
         s:=uniquedecomposition(f1);
        for j from 1 to m do
          if s[1]+a*s[2]=f[j][1] then
             P:=P*(s[1]+a*s[2])^(min(f[j][2],f2)); fi;
           if s[1]-a*s[2]=f[j][1] then
             P:=P*(s[1]-a*s[2])^(min(f[j][2],f2)); fi;
        od;
   od;P:=expand(P);rs:=[Re(P),Im(P)/sqrt(3)];
[rs,A*rs[1]+3*B*rs[2] mod 2*q,A*rs[2]-B*rs[1] mod 2*q]; end:
\end{verbatim}

In the case we have seen before, where $d=17$, if we take $a=1$,
$b=5$ and $c=29$, then we have $A=ac=29$, $B=bd=85$ and
$q=a^2+b^2=26$. The procedure above gives $r=-1$ and $s=7$ and
checks that $As+3Br\equiv 0$ (mod $2q$) and $Ar-Bs\equiv 0$ (mod
$2q$). These two conditions are enough to insure that the
expressions in (\ref{paramtwo}), all give integer values for the
coordinates of $\eta$ and $\zeta$. In the next procedure we use
the $r$ and $s$ determined earlier and construct an irreducible
regular octahedron in $\mathbb Z^3$, given a vector $(a,b,c)$ and
the possible values for $(m,n)$ in the decomposition of
$k^2=m^2-mn+n^2$. We are using simple formulae which can be
derived easily from Figure 1(b), thinking that the point $H$ is
the origin, $T[1]$, $T[2]$ and $T[3]$ are the points $M$, $N$ and
$L$. Then the other three vertices are simply given by the
addition of every two of these there vectors. Although there are
six possible equilateral triangles that one may start with in this
construction, one can see that in the end, essentially the same
regular octahedron (up to symmetries and translations) is
obtained.

\begin{verbatim}
> findpar:=proc(a,b,c,mm,nn)
local q,d,r,s,rs,A,B,k,mx,nx,my,ny,mz,nz,mu,nu,mv,nv,mw,nw,
u,v,w,x,y,z,T,R1,R2,DD,E,F;
  q:=a^2+b^2;k:=sqrt(mm^2-mm*nn+nn^2);d:=sqrt((a^2+b^2+c^2)/3);A:=a*c;B:=b*d;
    rs:=findgcd(A,-B,q);r:=rs[1][2]; s:=rs[1][1];
    mx:=-(d*b*(3*r+s)+a*c*(r-s))/(2*q);nx:=-(r*a*c+d*b*s)/q;
    my:=(d*a*(3*r+s)-b*c*(r-s))/(2*q);ny:=-(r*b*c-d*a*s)/q;mz:=(r-s)/2;nz:=r;
    mu:=nx;mv:=ny;mw:=nz;nu:=nx-mx;nv:=ny-my; nw:=nz-mz;
    u:=mu*mm-nu*nn;v:=mv*mm-nv*nn;w:=mw*mm-nw*nn;
    x:=mx*mm-nx*nn;y:=my*mm-ny*nn;z:=mz*mm-nz*nn;
    R1:=[(x+u-2*a*k)/3,(v+y-2*b*k)/3,(z+w-2*c*k)/3];
    R2:=[(x+u+2*a*k)/3,(v+y+2*b*k)/3,(z+w+2*c*k)/3];
 if R1[1]=floor(R1[1]) then T:=[[u,v,w],[x,y,z],R1]; else
                         T:=[[u,v,w],[x,y,z],R2];
 fi;
  DD:=[T[1][1]+T[2][1],T[1][2]+T[2][2],T[1][3]+T[2][3]];
  E:=[T[1][1]+T[3][1],T[1][2]+T[3][2],T[1][3]+T[3][3]];
  F:=[T[2][1]+T[3][1],T[2][2]+T[3][2],T[2][3]+T[3][3]];
 [T[1],T[2],T[3],DD,E,F];
end:
\end{verbatim}

Since for $k=2011$ we get $listofmn(k)=\{[880,2301]\}$, we checked
to see what regular octahedron is obtained for
$findpar(1,1,1,880,2301)$: [2301, -1421, -880], [880, -2301,
1421], [2401, 100, 1521], [3181, -3722, 541], [4702, -1321, 641],
and [3281, -2201, 2942]. Since the set $abcsol(2011)$ has 336
elements in it, there are essentially at most 337 irreducible
regular octahedra in $\mathbb Z^3$ with side lengths equal to
$2011\sqrt{2}$ as we pointed out in Proposition~\ref{numberofro}.
Next, we have a short function for subtracting two vectors $U$ and
$V$.

\begin{verbatim}
subtrv:=proc(U,V)
> local W;
  W[1]:=U[1]-V[1];W[2]:=U[2]-V[2];W[3]:=U[3]-V[3];[W[1],W[2],W[3]];
end:
\end{verbatim}
In order to compare various octahedrons it is easier if they are
all translated to the positive octant in such a way each plane of
coordinates contains at least a vertex of the octahedron. This is
accomplished with the next routine.

\begin{verbatim}
> tmttopqoctahedron:=proc(T)
local i,a,b,c,v,C;
 a:=min(T[1][1],T[2][1],T[3][1],T[4][1],T[5][1],T[6][1]);
 b:=min(T[1][2],T[2][2],T[3][2],T[4][2],T[5][2],T[6][2]);
 c:=min(T[1][3],T[2][3],T[3][3],T[4][3],T[5][3],T[6][3]);
   v:=[a,b,c];C:={subtrv(T[1],v),subtrv(T[2],v),subtrv(T[3],v),
      subtrv(T[4],v),subtrv(T[5],v),subtrv(T[6],v)};
end:
\end{verbatim}
So, for instance, the octahedron mentioned earlier of sides
lengths $2011\sqrt{2}$ becomes: [2401, 1521, 3822], [3822, 2401,
1521], [2301, 0, 1421], [1521, 3822, 2401], [0, 1421, 2301], and
[1421, 2301, 0]. The smallest cube $C_m=[0,m]^3$ containing an
octahedron positioned in the positive octant as specified earlier
is computed in the following procedure.

\begin{verbatim}
>mscofmoctahedron:=proc(Q)
 local a,b,c,T;
  T:=convert(Q,list);
   a:=max(T[1][1],T[2][1],T[3][1],T[4][1],T[5][1],T[6][1]);
   b:=max(T[1][2],T[2][2],T[3][2],T[4][2],T[5][2],T[6][2]);
   c:=max(T[1][3],T[2][3],T[3][3],T[4][3],T[5][3],T[6][3]);
 max(a,b,c);
 end:
\end{verbatim}

{\bf Step 3.} In our construction of these octahedrons we end up
with essentially the same octahedron if we proceed from a
different face of it. To eliminate the possibility of counting an
octahedron twice or more than one time, we would like to have a
way of distinguishing between octahedra and so an invariant to
translation and symmetries, like the side lengths, will be good.
Such an invariant is the set of $k$ values which are given by the
four pairs of opposite (parallel) faces of a regular octahedron.
We must have
$$\ell=d_1k_1\sqrt{2}=d_2k_2\sqrt{2}=d_3k_3\sqrt{2}=d_4k_4\sqrt{2}.$$
So, knowing $\ell$ and the $d_i$'s will give us $k_i$'s. The set
$\{k_1,k_2,k_3,k_4\}$ is clearly and invariant to translations and
the symmetries we have talked about at the beginning of the paper.
Hence, two octahedra with different sets of k-values will be
essentially different.  We determine these k-values from the first
three points given in the procedure $findpar$. The following
calculation finds $d/gcd(a,b,c)$, given $(a,b,c)$ such that
$a^2+b^2+c^2=3d^2$.

\begin{verbatim}
> unitvector:=proc(U)
local i,j,k,l,x;
   i:=U[1];j:=U[2];k:=U[3];
   l:=gcd(gcd(i,j),k); x:=(i^2+j^2+k^2)/(3*l^2);
sqrt(x);end:
\end{verbatim}

Then we need a routine to add three vectors.

\begin{verbatim}
> addvec:=proc(U,V,W)
local X;
X[1]:=U[1]+V[1]+W[1];X[2]:=U[2]+V[2]+W[2];X[3]:=U[3]+V[3]+W[3];
[X[1],X[2],X[3]];end:
\end{verbatim}

The next function is multiplying a scalar with a vector.

\begin{verbatim}
> multbyfactorv:=proc(v,k)
 local w; w[1]:=v[1]*k;w[2]:=v[2]*k;w[3]:=v[3]*k;
[w[1],w[2],w[3]];end:
\end{verbatim}

The Euclidean distance between to points is needed later.

\begin{verbatim}
> distance:=proc(A,B)
local C; C:=subtrv(A,B); sqrt(C[1]^2+C[2]^2+C[3]^2); end:
\end{verbatim}
As we said before, in the procedure that follows, $T$ is the list
of the first tthree vertices given by $findpar$.

\begin{verbatim}
> fourkvalues:=proc(T)
 local N1,N2,N3,N4,x,length;
 length:=distance(T[1],T[2])/sqrt(2);
   N1:=unitvector(addvec(T[1],T[2],T[3]));
   N2:=unitvector(addvec(T[1],T[2],multbyfactorv(T[3],-3)));
   N3:=unitvector(addvec(T[1],T[3],multbyfactorv(T[2],-3)));
   N4:=unitvector(addvec(T[3],T[2],multbyfactorv(T[1],-3)));
 {length/N1,length/N2,length/N3,length/N4};
 end:
\end{verbatim}

For the octahedron in Step 2, we have the set of $k$-values,
$\{1,2011\}$, as expected since $2011$ is a prime.

{\bf Step 4.}  In this step we calculate the orbit of an
octahedron $T$ within the reduced cube. The procedure $orbitbox1$
is taking care only of the eight possible symmetries.

\begin{verbatim}
> orbitbox1:=proc(T)
 local i,k,T1,a,b,c,x,T2,T3,T4,T5,T6,T7,T8,T9,T10,T11,T12,T13,T14,T15,T16,T17,T18,
 T19,T20,T21,T22,T23,T24,S,Q;
 Q:=convert(T,list);
   a:=max(Q[1][1],Q[2][1],Q[3][1],Q[4][1],Q[5][1],Q[6][1]);
   b:=max(Q[1][2],Q[2][2],Q[3][2],Q[4][2],Q[5][2],Q[6][2]);
   c:=max(Q[1][3],Q[2][3],Q[3][3],Q[4][3],Q[5][3],Q[6][3]);
   T1:=T;
   T2:={seq([Q[i][1],Q[i][2],c-Q[i][3]],i=1..6)};
   T3:={seq([Q[i][1],b-Q[i][2],Q[i][3]],i=1..6)};
   T4:={seq([a-Q[i][1],Q[i][2],Q[i][3]],i=1..6)};
   T5:={seq([a-Q[i][1],b-Q[i][2],Q[i][3]],i=1..6)};
   T6:={seq([a-Q[i][1],Q[i][2],c-Q[i][3]],i=1..6)};
   T7:={seq([Q[i][1],b-Q[i][2],c-Q[i][3]],i=1..6)};
   T8:={seq([a-Q[i][1],b-Q[i][2],c-Q[i][3]],i=1..6)};
   S:={T1,T2,T3,T4,T5,T6,T7,T8};
   S;
 end:
\end{verbatim}

The next procedure implements Theorem~\ref{coj} in our situation
where the objects of interest are octahedrons.

{\tiny
\begin{verbatim}
> orbitbox:=proc(T,p)
 local S,Q,TT,a,b,c,m,SS,i,d,u,v,w,y,mm,nn,x,k;
 Q:=convert(T,list);
 a:=max(Q[1][1],Q[2][1],Q[3][1],Q[4][1],Q[5][1],Q[6][1]);
 b:=max(Q[1][2],Q[2][2],Q[3][2],Q[4][2],Q[5][2],Q[6][2]);
 c:=max(Q[1][3],Q[2][3],Q[3][3],Q[4][3],Q[5][3],Q[6][3]);
 d:=max(a,b,c);u:=d-a;v:=d-b;w:=d-c;y:={u,v,w};k:=p-d;
    TT[1]:={seq([Q[i][3],Q[i][2],Q[i][1]],i=1..6)};
    TT[2]:={seq([Q[i][2],Q[i][3],Q[i][1]],i=1..6)};
    TT[3]:={seq([Q[i][1],Q[i][3],Q[i][2]],i=1..6)};
    TT[4]:={seq([Q[i][2],Q[i][1],Q[i][3]],i=1..6)};
    TT[5]:={seq([Q[i][3],Q[i][1],Q[i][2]],i=1..6)};
       S:=orbitbox1(T);
   for i from 1 to 5 do
       S:=S union orbitbox1(TT[i]);
   od;
 S:=convert(S,list);m:=nops(S);SS:={};
   for i from 1 to m do
     if S[i][1][1]<=a and S[i][2][1]<=a and S[i][3][1]<=a and S[i][4][1]<=a and S[i][5][1]<=a and S[i][6][1]<=a
        and S[i][1][2]<=b and S[i][2][2]<=b and S[i][3][2]<=b and S[i][4][2]<=b and S[i][5][2]<=b and S[i][6][2]<=b
        and S[i][1][3]<=c and S[i][2][3]<=c and S[i][3][3]<=c and S[i][4][3]<=c and S[i][5][3]<=c and S[i][6][3]<=c
        and S[i][1][1]>=0 and S[i][2][1]>=0 and S[i][3][1]>=0 and S[i][4][1]>=0 and S[i][5][1]>=0 and S[i][6][1]>=0
        and S[i][1][2]>=0 and S[i][2][2]>=0 and S[i][3][2]>=0 and S[i][4][2]>=0 and S[i][5][2]>=0 and S[i][6][2]>=0
        and S[i][1][3]>=0 and S[i][2][3]>=0 and S[i][3][3]>=0 and S[i][4][3]>=0 and S[i][5][3]>=0 and S[i][6][3]>=0
      then SS:=SS union {S[i]}; fi;
   od;
  mm:=nops(SS);
  nn:=(k+u+1)*(k+v+1)*(k+w+1);
    if nops(y)=3 then x:=6*mm*nn;fi;
    if nops(y)=2 then x:=3*mm*nn;fi;
    if nops(y)=1 then x:=mm*nn;fi;
 [x,SS,nops(SS),[u,v]]; end:
\end{verbatim}
}

The next two routines calculate the orbit of an octahedron within
its minimal cube $C_m$. This orbit has at most 48 elements and it
is needed in the process of elimination of octahedrons that have
already appeared in the construction. In comparison with the
previous orbit, it is bigger and invariant to all the symmetries.
\begin{verbatim}
> orbit1octahedron:=proc(T)
local i,k,T1,a,b,c,x,T2,T3,T4,T5,T6,T7,T8,T9,T10,T11,T12,T13,T14,
T15,T16,T17,T18,T19,T20,T21,T22,T23,T24,S,Q,d,a1,b1,c1;
 Q:=convert(T,list);
 d:=mscofmoctahedron(T);
 T1:=T;
 T2:={seq([Q[k][2],Q[k][3],Q[k][1]],k=1..6)};
 T3:={seq([Q[k][1],Q[k][3],Q[k][2]],k=1..6)};
 T4:={seq([Q[k][1],Q[k][2],d-Q[k][3]],k=1..6)};
 T5:={seq([Q[k][2],Q[k][3],d-Q[k][1]],k=1..6)};
 T6:={seq([Q[k][1],Q[k][3],d-Q[k][2]],k=1..6)};
 T7:={seq([Q[k][1],d-Q[k][2],Q[k][3]],k=1..6)};
 T8:={seq([Q[k][2],d-Q[k][3],Q[k][1]],k=1..6)};
 T9:={seq([Q[k][1],d-Q[k][3],Q[k][2]],k=1..6)};
 T10:={seq([d-Q[k][1],Q[k][2],Q[k][3]],k=1..6)};
 T11:={seq([d-Q[k][2],Q[k][3],Q[k][1]],k=1..6)};
 T12:={seq([d-Q[k][1],Q[k][3],Q[k][2]],k=1..6)};
 T13:={seq([Q[k][1],d-Q[k][2],d-Q[k][3]],k=1..6)};
 T14:={seq([Q[k][2],d-Q[k][3],d-Q[k][1]],k=1..6)};
 T15:={seq([Q[k][1],d-Q[k][3],d-Q[k][2]],k=1..6)};
 T16:={seq([d-Q[k][1],d-Q[k][2],Q[k][3]],k=1..6)};
 T17:={seq([d-Q[k][2],d-Q[k][3],Q[k][1]],k=1..6)};
 T18:={seq([d-Q[k][1],d-Q[k][3],Q[k][2]],k=1..6)};
 T19:={seq([d-Q[k][1],Q[k][2],d-Q[k][3]],k=1..6)};
 T20:={seq([d-Q[k][2],Q[k][3],d-Q[k][1]],k=1..6)};
 T21:={seq([d-Q[k][1],Q[k][3],d-Q[k][2]],k=1..6)};
 T22:={seq([d-Q[k][1],d-Q[k][2],d-Q[k][3]],k=1..6)};
 T23:={seq([d-Q[k][2],d-Q[k][3],d-Q[k][1]],k=1..6)};
 T24:={seq([d-Q[k][1],d-Q[k][3],d-Q[k][2]],k=1..6)};
 S:={T1,T2,T3,T4,T5,T6,T7,T8,T9,T10,T11,T12,T13,T14,T15,T16,
 T17,T18,T19,T20,T21,T22,T23,T24};
 S; end:
> orbitoctahedron:=proc(T)
 local S,Q,T1;
 Q:=convert(T,list);
 T1:={seq([Q[k][3],Q[k][2],Q[k][1]],k=1..6)};
 S:=orbit1octahedron(T) union orbit1octahedron(T1);
 S; end:
\end{verbatim}

\vspace{0.2in}

 {\bf Step 5.} At this point we are ready to build a list of minimal, irreducible
octahedrons. With each entry we keep, in order, $d$, the value $m$
which is  the size of the minimal cube ${\cal C}_m$, the six
vertices of the octahedron, the $k$-values and the corresponding
vector $(a,b,c)$. In the end, another list is created in a new
procedure, to get the list of all corresponding reducible cubes
needed in a calculation for a given dimension $n$. The
construction is a little complicated because there are octahedrons
for which all $k$-values are different of $1$.

\begin{verbatim}
> ExtendList:=proc(n,N,L::list,mm,nn,Orb::array)
local d,i,sol,nsol,nel,ttpcube,orb::array,
 LL::list,tnel,NL::list,cio,O,pO,k,kvalues,m,exception,Int,NN;
 NN:=floor((N+1)/2);
 orb:=array(1..2*NN+1);
 nel:=nops(L);
 LL:=L;
 k:=sqrt(mm^2-mm*nn+nn^2);
 m:=n*k;
 if m<=N then
  sol:=abcsol(n);nsol:=nops(sol);
  tnel:=nel;orb:=Orb;
    for i from 1 to nsol do
     O:=findpar(sol[i][1],sol[i][2],sol[i][3],mm,nn);pO:=tmttopqoctahedron(O);
         kvalues:=fourkvalues([O[1],O[2],O[3]]);
            exception:=evalb(1 in kvalues);
             if k=1 or exception=false then
                Int:=orbitoctahedron(pO) intersect orb[m];
                cio:=nops(Int);
                  if cio=0 then
                  orb[m]:=orb[m] union orbitbox(pO,N)[2];
                  d:=mscofmoctahedron(pO);                                                                                                                                                                                                                                                                                                                                                                                                                                                                                                                                                                                                                                                                               NL[tnel+1]:=[n,d,pO,kvalues,sol[i]];
              tnel:=tnel+1;
            fi;
       fi;
       od;
 LL:=[seq(L[i],i=1..nel),seq(NL[j],j=nel+1..tnel)];
 fi;LL,orb;
 end:
\end{verbatim}

 To generate all the other octahedrons we have to
magnify the irreducible ones.

\begin{verbatim}
> multbyfactoroctahedron:=proc(T,k)
local i,NT,Q;NT:={};
 Q:=convert(T,list);
 for i from 1 to 6 do
 NT:=NT union {multbyfactorv(Q[i],k)};
 od;NT; end:
\end{verbatim}

 The procedure $ExtendList$ is used recursively in the next loop
 to generate the list needed to calculate all octahedrons in ${\cal C}_N$ for a certain value $N$.

\begin{verbatim}
> ExtendListuptoN:=proc(N)
local i,j,k,l,kv,kvn,NN,L,Orb::array,mn,nmn,E,n,ii;
kv:=kvalues(N);kvn:=nops(kv);
   NN:=floor((N+1)/2);
   L:=[];Orb:=array(1..2*NN+1);
       for ii from 1 to 2*NN+1 do
         Orb[ii]:={};
       od;
     for i from 1 to kvn do
      k:=kv[i];mn:=listofmn(k);nmn:=nops(mn);
       for j from 1 to nmn do
        for l from 1 to NN+1 do
          n:=2*l-1;
           E:=ExtendList(n,N,L,mn[j][1],mn[j][2],Orb);
           L:=E[1];
          for ii from 1 to 2*NN+1 do
            Orb[ii]:=E[2][ii];
          od;
    od;
   od;
  od;
 L;end:
\end{verbatim}

As an example, for $n=20$, we obtain for $L:=ExtendListuptoN(20);$

[[1, 2, $\{[2, 1, 1], [1, 2, 1], [1, 1, 0], [0, 1, 1], [1, 1, 2],
[1, 0, 1]\}$, $\{1\}$, [1, 1, 1]],

[3, 4, $\{[4, 1, 0], [0, 3, 4], [3, 0, 4], [1, 4, 0], [4, 4, 3],
[0, 0, 1]\}$, $\{1, 3\}$, [1, 1, 5]],

[5, 10, $\{[8, 5, 7], [4, 10, 4], [1, 5, 8], [4, 0, 4], [7, 5, 0],
[0, 5, 1]\}$, $\{1\}$, [1, 5, 7]],

[7, 12, $\{[9, 0, 4], [12, 8, 9], [0, 4, 3], [3, 12, 8], [4, 3,
12], [8, 9, 0]\}$, $\{1, 7\}$, [1, 5, 11]],

[9, 16, $\{[3, 9, 0], [11, 7, 16], [14, 4, 4], [3, 0, 9], [0, 12,
12], [11, 16, 7]\}$, $\{1, 3\}$, [1, 11, 11]],

[11, 18, $\{[7, 18, 13], [11, 0, 1], [18, 7, 13], [15, 15, 0], [0,
11, 1], [3, 3, 14]\}$, $\{1\}$, [1, 1, 19]],

[13, 24, $\{[24, 15, 16], [16, 0, 9], [8, 24, 15], [15, 16, 0],
[9, 8, 24], [0, 9, 8]\}$, $\{1, 13\}$, [5, 11, 19]],

[13, 26, $\{[17, 13, 0], [24, 13, 17], [12, 26, 12], [7, 13, 24],
[12, 0, 12], [0, 13, 7]\}$, $\{1\}$, [7, 13, 17]],

[15, 28, $\{[15, 28, 13], [0, 19, 1], [20, 12, 0], [20, 9, 21],
[5, 0, 9], [0, 16, 22]\}$, $\{1, 3\}$, [5, 11, 23]],

[17, 24, $\{[0, 3, 4], [0, 20, 21], [13, 0, 24], [24, 4, 3], [11,
24, 0], [24, 21, 20]\}$, $\{1\}$, [1, 5, 29]],

[17, 34,$\{[23, 17, 0], [15, 34, 15], [15, 0, 15], [0, 17, 7],
[30, 17, 23], [7, 17, 30]\}$, $\{1\}$, [7, 17, 23]],

[19, 36, $\{[0, 12, 12], [34, 24, 24], [11, 17, 36], [11, 36, 17],
[23, 0, 19], [23, 19, 0]\}$, $\{1\}$, [5, 23, 23]],

[19, 30, $\{[0, 21, 25], [9, 25, 0], [21, 5, 30], [30, 9, 5], [25,
30, 21], [5, 0, 9]\}$, $\{1, 19\}$, [1, 11, 31]]]

\n Let us observe that although there are four primitive solutions
for

$$abcsol(19)=[[5, 23, 23], [11, 11, 29], [13, 17, 25], [1, 11,
31]]$$

\n we get only two essentially  different octahedrons of side
lengths $19\sqrt{2}$. We also need to take into account the
reducible octahedrons.

\begin{verbatim}
> L:=ExtendListuptoN(100):
 ExtendListuptoNmultiples:=proc(N,L)
 local i,j,x,lc,m,mm,d,dd,C,CC,LL;
 m:=nops(L);i:=1;LL:={};
 while i<=m  do
   d:=L[i][2];
      if d<=N then
         mm:=floor(N/d);C:=L[i][3];j:=2;
          while j<=mm  do
               CC:=multbyfactoroctahedron(C,j);
               dd:=d*j;lc:=nops(LL);
               LL:=LL union {[L[i][1]*j,dd,CC,L[i][4]]};
              j:=j+1;
          od;
       fi;
     i:=i+1;
  od;
 convert(LL,list);
 end:

> LL:=ExtendListuptoNmultiples(20,L):
\end{verbatim}

The program is taking the previous list and inflates only the
octahedrons needed:

LL := [[3, 6, $\{[3, 3, 6], [3, 3, 0], [3, 6, 3], [0, 3, 3], [6,
3, 3], [3, 0, 3]\}, \{1\}]$,

[2, 4, $\{[2, 2, 4], [2, 2, 0], [2, 4, 2], [0, 2, 2], [4, 2, 2],
[2, 0, 2]\}, \{1\}$],

[7, 14, $\{[7, 7, 14], [7, 7, 0], [7, 14, 7], [0, 7, 7], [14, 7,
7], [7, 0, 7]\}, \{1\}$],

[6, 12, $\{[6, 6, 12], [6, 6, 0], [6, 12, 6], [0, 6, 6], [12, 6,
6], [6, 0, 6]\}, \{1\}$],

[5, 10, $\{[5, 10, 5], [0, 5, 5], [5, 5, 10], [5, 5, 0], [10, 5,
5], [5, 0, 5]\}, \{1\}$],

[4, 8, $\{[4, 0, 4], [4, 4, 0], [4, 4, 8], [4, 8, 4], [0, 4, 4],
[8, 4, 4]\}, \{1\}$],

[10, 20, $\{[10, 10, 20], [10, 10, 0], [10, 20, 10], [0, 10, 10],
[10, 0, 10], [20, 10, 10]\}, \{1\}$],

[9, 18, $\{[9, 9, 18], [9, 9, 0], [9, 18, 9], [0, 9, 9], [18, 9,
9], [9, 0, 9]\}, \{1\}$],

[8, 16, $\{[8, 8, 16], [8, 8, 0], [8, 16, 8], [0, 8, 8], [16, 8,
8], [8, 0, 8]\}, \{1\}$],

[15, 20, $\{[0, 15, 20], [0, 0, 5], [15, 0, 20], [5, 20, 0], [20,
20, 15], [20, 5, 0]\}, \{1, 3\}$],

[12, 16, $\{[0, 12, 16], [4, 16, 0], [16, 4, 0], [0, 0, 4], [12,
0, 16], [16, 16, 12]\}, \{1, 3\}$],

[9, 12, $\{[0, 9, 12], [0, 0, 3], [9, 0, 12], [12, 12, 9], [3, 12,
0], [12, 3, 0]\}, \{1, 3\}$],

[6, 8, $\{[2, 8, 0], [8, 8, 6], [8, 2, 0], [0, 6, 8], [6, 0, 8],
[0, 0, 2]\}, \{1, 3\}$],

[10, 20, $\{[8, 0, 8],[14, 10, 0], [8, 20, 8], [0, 10, 2], [16,
10, 14], [2, 10, 16]\}, \{1\}$]]

\vspace{0.2in}

 {\bf Step 6}. Finally we are adding up the contribution
of each octahedron located in the union of the two lists created
earlier.

\begin{verbatim}
> addupnew:=proc(N,L,LL)
local i,j,k,nc,mm,m,d,C,CC,x,dd,nt; nc:=0; m:=nops(L); i:=1;
   while i<=m  do
   d:=L[i][2];
            if d<=N then
                x:=orbitbox(L[i][3],N)[1];
                nc:=nc+x;
            fi;  i:=i+1;
   od;
 m:=nops(LL);
  i:=1;
    while i<=m  do
    d:=LL[i][2];
      if d<=N then
              x:=orbitbox(LL[i][3],N)[1];
              nc:=nc+x;
      fi;  i:=i+1;
   od;
 nc;end:
\end{verbatim}

The next command produces the sequence we are looking for.
\begin{verbatim}
 >NO:=[seq([k,addupnew(k,L,LL)],k=1..100)];
\end{verbatim}

So, the first one hundred terms of A178797 are:

NO := [[1, 0], [2, 1], [3, 8], [4, 32], [5, 104], [6, 261], [7,
544], [8, 1000], [9, 1696], [10, 2759], [11, 4296], [12, 6434],
[13, 9352], [14, 13243], [15, 18304], [16, 24774], [17, 32960],
[18, 43223], [19, 55976], [20, 71752], [21, 90936], [22, 113973],
[23, 141312], [24, 173436], [25, 210960], [26, 254587], [27,
305000], [28, 364406], [29, 432824], [30, 511421], [31, 600992],
[32, 702556],[33, 817200], [34, 946131], [35, 1090392], [36,
1251238], [37, 1430072], [38, 1629391], [39, 1850064], [40,
2094276], [41, 2363616], [42, 2659813], [43, 2984600], [44,
3341660], [45, 3731720], [46, 4156689], [47, 4618480], [48,
5119292], [49, 5661600], [50, 6248705], [51, 6882808], [52,
7568126], [53, 8306520], [54, 9104339], [55, 9962320], [56,
10888762], [57, 11882896], [58, 12949661], [59, 14090952], [60,
15311286], [61, 16613736], [62, 18001975], [63, 19479680], [64,
21052826], [65, 22724576], [66, 24500175], [67, 26383240], [68,
28387456], [69, 30510616], [70, 32758963], [71, 35136544], [72,
37656214], [73, 40317328], [74, 43125329], [75, 46085496], [76,
49207224], [77, 52493112], [78, 55954267], [79, 59592272], [80,
63415296], [81, 67428832], [82, 71642127], [83, 76059704], [84,
80701546], [85, 85565064], [86, 90662451], [87, 95997360], [88,
101592122], [89, 107443264], [90, 113561009], [91, 119951832],
[92, 126644136], [93, 133629672], [94, 140916757], [95,
148513712], [96, 156444624], [97, 164706400], [98, 173308509],
[99, 182260568], [100, 191575248]]


\begin{thebibliography}{99}
\bibitem{a} N. C. Ankeny, {\it Sums of Three Squares}, {\em Proceedings of
AMS}, vol. 8, No. 2, pp 316-319.

\bibitem{rceji} R.~Chandler and E.~J.~Ionascu,  {\it A characterization of all equilateral triangles in $\Bbb Z^3$},
{\em Integers}, Art. A19 of Vol. {\bf 8} 2008.

\bibitem{ch} S. Cooper and M. Hirschhorn, On the number of primitive
representations of integers as a sum of squares, {\em Ramanujan J.}
(2007) 13, pp. 7-25.

\bibitem{cox} D. A. Cox, Primes of the Form $x^2 + ny^2$: Fermat, Class Field Theory, and Complex
Multiplication, Wiley-Interscience, 1997

\bibitem{g} R. Guy, {\em Unsolved Problems in Number Theory}, Springer-Verlag, 2004
\bibitem{gr} E. Grosswald, {\it Representations of integers as sums of squares}, Springer Verlag, New York, 1985.
\bibitem{s} I. J. Schoenberg, {\it Regular Simplices and Quadratic Forms}, J. London Math. Soc. 12 (1937) 48-55.

\bibitem{hs} M. D. Hirschhorn and J. A. Sellers, On representations of numbers as a sum of three squares,
{\em Discrete Mathematics}, 199 (1999), pp. 85-101.

\bibitem{eji} E. J. Ionascu, {\em A parametrization of equilateral triangles having integer
coordinates}, Journal of Integer Sequences, Vol. {\bf 10}, 09.6.7.
(2007)

\bibitem{ejic} E. J. Ionascu, {\em Counting all equilateral triangles in $\{0,1,2,...,n\}^3$}, Acta
Math. Univ. Comenianae, vol. LXXVII, {\bf 1}(2008), pp.129-140.

\bibitem{ejicregt} E. J. Ionascu, {\em Counting all regular tetrahedra in $\{0,1,2,...,n\}^3$}, submitted for
publication, arXiv:0912.1062v1

\bibitem{ejio} E. J. Ionascu and R. Obando {\em Counting all cubes in $\{0,1,2,...,n\}^3$}, submitted for
publication, arXiv:1003.4569

\bibitem{ejirt} E. J. Ionascu, {\em A characterization of regular tetrahedra in $\Bbb Z^3$},
J. Number Theory, 129(2009), 1066-1074.

\bibitem{ejiam} E. J. Ionascu and A. Markov, {\em Platonic solids in $\mathbb
Z^3$}, to appear in J. Number Theory, arXiv:0910.1722v1


\bibitem{il} I. Larrosa, {\em Solution to Problem 8,
http://faculty.missouristate.edu/l/lesreid/POW08$\_$03.html}

\bibitem{r} K. Rosen, {\em Elementary Number Theory}, Fifth Edition, Addison Wesley,
2004.

\bibitem{OL} Neil J.~A. Sloane, \newblock {\em The On-Line Encyclopedia of
Integer Sequences}, \newblock 2005, \newblock published
electronically at http://www.research.att.com/$\sim$njas/sequences/.
\end{thebibliography}
\end{document}